\documentclass[a4paper,10.2pt]{article}

\usepackage[latin1]{inputenc}
\usepackage[T1]{fontenc}
    \usepackage[english]{babel}
    \usepackage{amsmath, amssymb}
\usepackage{geometry,amsmath,amssymb,enumerate,bbm,amsthm}
\usepackage{amsfonts}
\usepackage{amsmath,amsthm,indentfirst}
\usepackage{amssymb}
\usepackage{pstricks}
\usepackage{amsthm}
\usepackage[matrix,arrow]{xy}
\usepackage{amssymb}
\usepackage{amsmath}
\usepackage{vmargin}
\usepackage{dsfont}
\usepackage{mathrsfs}
\usepackage{fancyhdr}

\usepackage[pdftex,bookmarks,colorlinks,breaklinks]{hyperref}  % PDF hyperlinks, with coloured links

\usepackage{memhfixc}  % remove conflict between the memoir class & hyperref
\usepackage{pdfsync}  % enable tex source and pdf output syncronicity

\newcommand{\eps}{\epsilon}

\newcommand{\ep}{\epsilon}
\everymath{\displaystyle}

\newcommand{\RR}{\mathbb{R}}

\def\qed{\hbox{${\vcenter{\vbox{
  \hrule height 0.4pt\hbox{\vrule width 0.4pt height 6pt
  \kern5pt\vrule width 0.4pt}\hrule height 0.4pt}}}$}}

\newtheorem{theo}{Theorem}
\newtheorem*{theo'}{Theorem 8 (Resnick)} 
\newtheorem*{theo''}{Theorem 9 (Abidi and Hmidi)}
\newtheorem{prop}[theo]{Proposition}
\newtheorem{lemm}[theo]{Lemma}

\newtheorem{remark}[theo]{Remark}

\newtheorem{defi}[theo]{Definition}

\date{}
\author{Omar Lazar}
\title{Global and local existence for the dissipative critical SQG equation with small oscillations } 
\begin{document}
\maketitle
\bibliographystyle{plain}
\noindent{\bf Abstract:} This article is devoted to the study of the critical dissipative surface quasi-geostrophic $(SQG)$ equation in $\mathbb{R}^2$. For any initial data $\theta_{0}$
belonging to the space $\Lambda^{s} ( H^{s}_{uloc}(\mathbb{R}^2)) \cap L^\infty(\mathbb{R}^2)$, we show that  the critical (SQG) equation has at least one global weak solution in time for all $1/4\leq s \leq 1/2$ and at least one local weak solution in time for all $0<s<1/4$.  The proof  for the global existence is based on a new energy inequality which improves the one obtain in \cite{Laz} whereas the local existence uses more refined energy estimates based on Besov space techniques. 

\section{Introduction}

In this article, we study the following two dimensional dissipative surface quasi-geostrophic equation:
\begin{equation}%\label{qg}
\ (SQG)_{\alpha}: \\\left\{
\aligned
&\partial_{t}\theta(x,t)+u.\nabla\theta+ \nu\Lambda^{\alpha}\theta = 0,
\\ \nonumber
& u(\theta)= \mathcal{R}^\perp \theta,
\\ \nonumber
& \theta(0,x)=\theta_{0}(x),
\endaligned
\right.
\end{equation}
where $\nu>0$ is the viscosity, which we will assume to be equal to 1 without loss of generality, and  
\[\Lambda^{\alpha}\theta\equiv (-\Delta)^{\alpha/2} \theta=C_{\alpha}P.V.\int_{\mathbb{R}^{2}}{\frac{\theta(x)-\theta(x-y)}{|y|^{2+\alpha}}dy},\]
where $C_{\alpha}$ is a positive constant and $\theta\,:\,\mathbb R_{+}\times \mathbb R^2\to \mathbb R$ is a scalar function which modelises the potential temperature of the fluid (typically ocean).  Here $\alpha \in (0,2]$ is a fixed parameter and the velocity
$u=(u_1,u_2)$ is divergence-free and determined by the Riesz
transforms of the potential temperature $\theta$ via the formula:
$$
u=(-{\mathcal R}_2\theta,{\mathcal R}_1\theta)=(-\partial_{x_2}
(-\Delta)^{-1/2}\theta,
\partial_{x_1}(-\Delta)^{-1/2}\theta).
$$
It is well known that this equation is classified into 3 cases depending on the value of the power of the fractional Laplacian. Namely, the cases $0<\alpha<1$, $\alpha=1$ and $1<\alpha<2$  are respectively called super-critical, critical, and sub-critical case. Before recalling some well known results, let us mention that this equation was introduced in 1994 by Constantin, Majda and Tabak \cite{CMT} in order to get a better understanding of the 3D Euler equation. More precisely, the authors in \cite{CMT} have pointed out that, in the inviscid case (i.e $\nu=0$), the (SQG) equation written in terms of the gradient $\nabla^{\perp}\theta$ turns out to share the same properties with the 3D Euler equation written in terms of vorticity (i.e the curl of the velocity field).  It was conjectured in \cite{CMT} that the gradient of the active scalar has a fast growth when the geometry of the level sets contain a hyperbolic saddle. In \cite{C}, C\'ordoba showed that a simple hyperbolic saddle breakdown cannot occur in finite time.
Beside being mathematically relevant, the inviscid (SQG) equation appears also in some physical models such as the study of strongly rotating fluids.  In fact, the inviscid (SQG) equation comes from more general quasi-geostrophic models of atmospheric and ocean fluid (see e.g \cite{MB}, \cite{Ped}). Even in the viscous case, and more precisely in the critical case $\alpha=1$, the equation is physically pertinent. In fact, the term $(-\Delta)^{1/2} \theta$ modelises the so-called Eckman pumping which, roughly speaking, drives geostrophic flows. The existence of global weak $L^{2}$ solutions for the inviscid and dissipative $(SQG)_{\alpha}$ equation goes back to Resnick \cite{Res}. Subsequently, Marchand \cite{Mar} has extended Resnick's result to the class of initial data belonging to $H^{-1/2}$ or $L^{p}$ with $p>4/3$. \\
The sub-critical case ($1<\alpha<2$) is well understood, for instance one can cite the works by Resnick \cite{Res} and  latter Ju \cite{Ju} in which they  showed that for all $\theta_{0}\in H^{s}$, with $s>2-2\alpha$, and $s\geq 0$ there exists a unique global smooth solution. \\
 In the critical case, namely $\alpha=1$, Constantin, C\'ordoba and Wu \cite{CCW}  showed that there is a unique global solution when $\theta_0$ is in the critical space  $H^{1}$  under a smallness assumption on $\Vert \theta_0 \Vert_{\infty}$. The problem consisting in showing the regularity of those weak $L^{2}$ solutions has been solved by Caffarelli and Vasseur \cite{CV}  using De Giorgi iteration method. By using a completely different method, Kiselev, Nazarov and Volberg \cite{KNV} showed that all periodic smooth initial data give rise to a unique smooth solution. Their proof is based on a new non local maximum principle verified by the gradient of the solution which is shown to be bounded from above by the derivative of the modulus of continuity at time 0.  Both methods have been succesfully applied in several works. Abidi and Hmidi \cite{Hm} showed that for all $\theta_{0} \in \dot B^{0}_{\infty, 1}$ there exist a unique global solution by using a Langragian approach and the method of \cite{KNV} for showing that local solutions are global. See also the work by Dong, Li \cite{DL1} and Dong, Du \cite{DL2} where same type of results have been obtained. \\
 As for the supercritical case ($0<\alpha<1$), only partial results are known. We can cite the result of C\'ordoba and C\'ordoba \cite{CC} where global solutions were shown to exist for all small $H^m$ initial data with $m<2$. This latter result was improved by Chae and Lee \cite{CL} where global existence and uniqueness are shown for data in critical Besov space $B^{2-\alpha}_ {2,1}$ under a smallness assumption of the $\dot B^{2-\alpha}_ {2,1}$ norm. In \cite{Ju}, Ju established the global regularity for small $H^{s}$ data  $s\geq2-\alpha$ and therefore improved the two previous results. We can also cite the work of Chen, Miao and Zhang \cite{CMZ} in which they proved the global well-posedness for small initial in critical Besov spaces and local well-posedness for large data.  Constantin and Wu  studied the H\"older regularity propagation in \cite{CW} where they proved that all weak $L^2$ solutions  which have the property to be 
 $\delta$-H\"older continuous (with $\delta>1$) on the time interval $[t_1,t]$,  are a classical solutions on $(t_1,t]$. The result obtained in \cite{CW} is sharp.  Indeed, Sylvestre, Vicol and Zlato\v{s}  in \cite{SVZ}   proved the existence of a solution with time-independent velocity $u \in \mathcal{C}^{\delta}$, with $\delta \in (0,1)$ and small $C^{2}$ initial data  that becomes discontinuous in finite time. It is an outstanding open problem to prove that any bounded weak solutions become $\mathcal{C}^{1-\alpha}$ after a short time or if they blow-up in finite time, uniqueness of these solutions is also open. \\

In this paper, we present some existence results for the critical dissipative (SQG) equation with data in the non homogeneous Morrey-Campanato  spaces $M_{p,q} (\mathbb{R}^{d})$ ($1\leq q \leq p \leq \infty$) defined  by:
$$M_{p,q} (\mathbb{R}^{d})=\left\{ f \in L^{q}_{loc} (\mathbb{R}^{2}), \ \sup_{x\in \mathbb{R}^{n}}\sup_{0<R<1} R^{d(\frac{1}{p}-\frac{1}{q})}\Vert f \Vert_{L^{q}(B_{R}(x))} < \infty\right \} $$
 In the same spirit as the work by Lemari\'e-Rieusset for the Navier-Stokes equation cases, we study the critical (SQG) equation with data in the space $L^{2}_{uloc}(\mathbb{R}^{2})$ (see \cite{PGLR}, \cite{PGLR2}). Other authors have obtained existence result for the Navier-Stokes equation with data in homogeneous Morrey-Campanato spaces (see e.g. Basson \cite{Basson}, Federbush \cite{Fed}, or  Taylor \cite{Tayl}). Unfortunately, because of the Riesz operator, this space is not well adapted to the critical (SQG) equation. This is due to the fact that the Riesz transforms are not defined in the space $L^{2}_{uloc}(\mathbb{R}^{2})$ they are not even defined in the subspace $L^{2}_{loc}(\mathbb R^2)$. In fact, the problem comes from the lack of integrablity at infinity of the Riesz kernel. In \cite{Laz}, we have overcome the difficulty by putting some integrability on the solution, that is to say, we have worked with $\Lambda^{s} \theta$ instead of $\theta$, where $0<s<1$. In this new setting, we were able to work in a space close to the space $L^{2}_{uloc} (\mathbb R^2)$, namely, the space $\Lambda^{s} ( H^{s}_{uloc}(\mathbb R^2))$ (where $f\in \Lambda^{s}\left.( H^{s}_{uloc}(\mathbb R^2)\right.)$ if there exists $g\in   H^{s}_{uloc}(\mathbb R^2)$ s.t $f=\Lambda^{s} g$).  In the case $1/2< s < 1$ we have obtained global existence of  weak solutions in \cite{Laz}.
 In this article, by slightly improving the energy inequality obtained in \cite{Laz}, we are able to show that the solutions are in fact global for all $1/4\leq s \leq 1/2$. For the case $0<s<1/4$, we only have a local existence result. The case $s=1/4$ turns out to be critical in the sense that below this regularity exponent, we have a lack of regularity to control the $L^{\infty}$ norm of the Riesz transforms (the problem comes from the control of the low frequencies). It is worth mentioning that the smaller the value of $s$, the closer the singular case  i.e $L^{2}_{uloc} (\mathbb R^2)$ data, and the  partial result obtained for small $s$ is in part due to this fact.  \\

 The paper is organized as follows. In the first part we recall some tools we will use in the proof of our main results. We also briefly recall the notations,  and we present the equation with truncated and regularized initial data. In the second part, we prove an energy inequality available in the case $1/4\leq s \leq 1/2$ which allows us to get a global control of the solutions. The third part is devoted to the proof of an energy inequality available in the case $0< s <1/4$, this inequality leads to a local existence theorem. The passage to the limit is omitted here since it is the same as the one in \cite{Laz}.  \\
 
 Throughout this article, $C$ stands for any controlled and positive constant, which therefore could be different from line to line. We also denote $A \lesssim B$ if $A$ is less to $B$ up to a positive multiplicative constant which can be different from line to line as well. Those constants depend only on some controlled norms. We denote by $\mathcal{D}(\mathbb R^{2})$ the space of smooth functions in $\mathbb R^{2}$ that are compactly supported. As usually, we denote by $L^{p}$ the usual Lebesgue spaces, $H^{s}$  the classical Sobolev spaces and $\dot H^{s}$ the homogeneous ones.  We shall also use the shorter notation $L^{p}X$ for $L^{p}([0,T],X)$ where $X$ is a Lebesgue or Sobolev space.   \\

 Our  main result reads as follows.
 \begin{theo} \label{tp} Let us denote $X_{T}$ and  $X^{s}_{T}$ the spaces
$$
X_{T}\equiv L^\infty ([0,T], L^{2}_{uloc})\cap (L^{2}_{t}([0,T], \dot H^{1/2}))_{uloc},
$$
and 
$$
X^{s}_{T} \equiv L^\infty ([0,T], \dot H^{s}_{uloc})\cap (L^{2}_{t}([0,T], \dot H^{s+1/2}))_{uloc}.
$$ \\
Assume that  $\theta_0= \Lambda^{s} w_{0} \in \Lambda^{s}(H^{s}_{uloc}) \cap L^{\infty},$ then: \\

$\bullet$  If $1/4\leq s\leq1/2$,  the critical (SQG) has at least one global weak solution $\theta$ which satisfies $\theta \in X_T$ and $w \in X^{s}_{T}$ for all $T<\infty$. Futhermore, we have the following control

 $$
\Vert w(x,t) \Vert^{2}_{{\dot H^{s}_{uloc}(\mathbb{R}^2)}} \leq c \ e^{CT},
$$

$\bullet$  If $0<s<1/4$,  the critical (SQG) has at least one local weak solution $\theta$ so that for all $$T<T^{*}\equiv {\frac{C(\Vert \theta_0 \Vert_{\infty})} {1+\Vert w_{0} \Vert^{2}_{\dot H^{s}_{uloc}}}}$$ 
we have
$$\theta \in X_{T^*} \ \text{and} \ w \in X^{s}_{T^*}.$$

\noindent  Moreover, for all $T\leq T^*$, the solution $w$ satisfies the following energy inequality:
$$
\Vert w(x,T) \Vert^{2}_{H^{s}_{uloc}} \leq \Vert w_{0} \Vert^{2}_{ H^{s}_{uloc}} +   C \int_{0}^{T} \left(\Vert w(x,s) \Vert^{2}_{\dot H^{s}_{uloc}}+\Vert w(x,s) \Vert^{4}_{{\dot H^s}_{uloc}}  +  C \right) \ ds,
$$

\noindent  where $C$ is a positive constant depending only on  $\displaystyle\Vert \theta_{0} \Vert_{L^{\infty}(\mathbb{R}^2)}$ and \ $\Vert w_{0} \Vert_{ H^{s}_{uloc}(\mathbb{R}^2)}$.
\end{theo}
\begin{remark}
Previous results obtained in \cite{Laz} and the global result  of the above theorem imply that we have global existence for $1/4 \leq s \leq 1$.
\end{remark}

\begin{remark}
In fact, for the case $0<s< 1/4$, we prove a more general inequality, namely, for all $t\leq T^*$ and for all $K \in (0,2)$,  solution $w$ satisfies the following energy inequality
$$
\Vert w(x,T) \Vert^{2}_{ H^{s}_{uloc}} \leq \Vert w_{0} \Vert^{2}_{ H^{s}_{uloc}} +   C \int_{0}^{T} \left(\Vert w \Vert^{2}_{\dot H^{s}_{uloc}}+\Vert w \Vert^{2(\frac{3-K}{2-K})}_{{\dot H^s}_{uloc}}  +  C \right) \ ds.
$$
The best power that we can obtain in the  above inequality is close to 3 (roughly speaking, it corresponds to $K \rightarrow 0$) and therefore we only have a local control of the solutions. In order to make the statement and the time existence of the solutions simpler, we have considered the case $K=1$ which gives a power 4. 
\end{remark}
 
 \section{The $L^{p}_{uloc} (\mathbb R^2 )$ and  $H^{s}_{uloc}(\mathbb R^2)$ spaces}
 In this section, we recall the definition of the $L^{p}_{uloc} (\mathbb R^2 )$ and $H^{s}_{uloc} (\mathbb R^2)$ spaces. To do so, we need to introduce the set of translations of a given test function.
 \begin{defi}
Let us fix a positive test function $\phi_{0}$ such that $ \phi_{0} \in \mathcal{D} (\mathbb R^2) $ and
\begin{equation}
\left \{
\aligned
&\phi_{0}(x)= 1 &   \mathrm{if} \ \vert x\vert \leq 2,
\\ \nonumber
& \phi_{0}(x)= 0 &   \mathrm{if} \  \vert x \vert \geq 3.
\endaligned
\right.
\end{equation}
\end{defi}
\noindent We define the set of translations of the function $\phi_{0}$ as $B_{\phi_{0}}\equiv \{\phi_{0}(x-k), k\in \mathbb Z^2 \}.$  We are now ready to define both $L^{p}_{uloc} (\mathbb R^2 )$ and $H^{s}_{uloc} (\mathbb R^2)$ spaces.
\begin{defi}
Let $1\leq p \leq \infty$ then  $f \in L^{p}_{uloc}(\mathbb R^2)$   if and only if  $f \in L^{p}_{loc}(\mathbb R^2)$ and the following norm is finite
 $$
\Vert f \Vert_{L^{p}_{uloc}(\mathbb R^2)}= \sup_{\phi \in {B_{\phi_0}}} \Vert \phi f \Vert_{L^{p}(\mathbb R^2)}.
$$
\end{defi}

\noindent We will also use the following useful equivalent norms
$$\Vert f \Vert_{L^{p}_{uloc}(\mathbb R^2)} \approx \sup_{k \in \mathbb{Z}^2 } \left(\int_{k+[0,1]^2} \vert f(x) \vert^p \ dx \right)^{1/p} \approx \sup_{k \in \mathbb{Z}^2 } \Vert \phi_{0} (x-k) f \Vert_{L^{p}(\mathbb R^2)}.  $$
Let us recall the definition of the $H^{s}_{uloc}(\mathbb R^2)$ spaces with $0<s<1$.
 \begin{defi}
Let $\phi_{0}$ be a positive test function chosen as in definition 4. We say that $f \in H^{s}_{uloc}(\mathbb R^2)$   if and only if  $f \in H^{s}_{loc}(\mathbb R^2)$ and  the following norm is finite
\begin{center}
$\displaystyle \Vert f \Vert^{2}_{{H^{s}_{uloc}} (\mathbb R^2)}  = \sup_{\phi\in B_{\phi_0}  } \Vert \phi f \Vert_{H^{s}}.$
\end{center}
\end{defi}
\noindent We will also use the following equivalent  norms, let us set 
\begin{equation} \label{norm}
A_{\phi} f\equiv\int \frac{ \vert \phi f \vert^2}{2} + \frac{ \vert \Lambda^{s}(\phi f) \vert^{2}}{2} \ dx,
\end{equation}
then,
\begin{center}
$\displaystyle \Vert f \Vert^{2}_{{H^{s}_{uloc}} (\mathbb R^2)} = \sup_{\phi\in B_{\phi_0}  }  A_{\phi} f,$
\end{center}
and,
$$
\Vert f \Vert^{2}_{\dot H^{s}_{uloc}}= \sup_{k \in \mathbb{Z}^{2}} \int \vert \Lambda^{s} (\phi_{0}(x-k) f(x,t)) \vert^{2} \ dx.
$$
\begin{remark} It is important to note that the norms do not depend on the choice of the test functions. In \cite{Laz}, we have seen that the $H^{s}_{uloc} (\mathbb R^{2})$ can be considered with $\phi$ inside or outside the fractional derivative since the norms are equivalents. Therefore, it does not matter whether the function $\phi$ is inside or outside the brackets in our computations.
\end{remark}
The spaces  $(L^{2}_{T}\dot H^{s})_{uloc}$ and  $L^{\infty}_{T} \dot H^{s}_{uloc}$ with $0<s<1$ and $1\leq p < \infty$ will be used throughout the paper. These spaces are endowed with the following norms
\begin{eqnarray*}
\Vert w \Vert^{2}_{ (L^{2}_{T} \dot H^{s})_{uloc}}&=& \sup_{\phi \in B_{\phi_0}} \int_{0}^{T} \int \phi \vert \Lambda^{s} w(x,s) \vert^2 \ dx \ ds < \infty, \\
\Vert w \Vert^{2}_{L^\infty_{T} \dot H^{s}_{uloc}} &=& \sup_{t\in[0,T]} \sup_{\phi \in B_{\phi_0}} \int \phi \vert \Lambda^{s} w(x,t) \vert^{2} \ dx  < \infty.
\end{eqnarray*}

In the next section, we recall some classical tools from the so-called Littlewood-Paley theory (see e.g \cite{MC}, or \cite{PGLR}).
\section{Besov spaces and Bernstein inequalities.}
In order to define the Besov spaces, we need to recall the definition of the dyadic blocks.
\begin{defi}
Let $\phi \in \mathcal{D}(\mathbb R^2)$ be a non negative function such that $\phi(x)=1$ if $\vert x \vert \leq 1/2$ and $0$ if $\vert x \vert \geq 1$. We also define $\psi \in \mathcal{D}(\mathbb R^2)$ by $\psi(x)=\phi(x/2)-\phi(x)$ supported on a corona. Then, we define the Fourier multiplier $S_j$ and $\Delta_j$ by
$$
\widehat{S_{j} f}(\xi)=\phi(\frac{\xi}{2^{j}}) \widehat{ f}(\xi) \ \ \text{and} \ \  \widehat{\Delta_{j} f}(\xi)=\psi(\frac{\xi}{2^{j}}) \widehat{ f} (\xi).
$$
\end{defi}
\hspace{-0,5cm}From these operators we deduce the Littlewood decomposition of a distribution  $f\in \mathcal{S}'$, that is, for all $N\in \mathbb Z$, we have
$$
f=S_N f+ \sum_{j\geq N} \Delta_j f \ \text{in} \ \mathcal{S}'(\mathbb{R}^{2}).
$$
If moreover,
$$
S_N  f \xrightarrow[N\to-\infty]{ } 0 \ \ \text{in} \ \ \mathcal{S}'(\mathbb{R}^{2}),
$$
we obtain the homogeneous decomposition of $f \in \mathcal{S}'(\mathbb{R}^{2}) $ (modulo polynomials)  :
$$
f= \sum_{j\in \mathbb{Z} } \Delta_j f \hspace{0,6cm} \text{in} \hspace{0,2cm} \ \mathcal{S}'(\mathbb{R}^{2}).
$$
The inhomogeneous Besov spaces are defined as follow
\begin{defi}
For $s\in \mathbb R$, $(p, q) \in [1, \infty]^2$ and $N \in \mathbb Z$, a distribution $f \in \mathcal{S}'(\mathbb R^2)$ belongs to the inhomogeneous Besov space $B^{s}_{p,q}$ if and only if
$$\Vert f \Vert_{B^{s}_{p,q}}=\Vert \Delta_N f \Vert_{p} + \left( \sum_{j\geq N} 2^{jqs} \Vert \Delta_j f \Vert^{q}_{p} \right)^{1/q} < \infty $$
\end{defi}
\noindent We also recall the definition of the homogeneous Besov spaces  which are defined only modulo polynomials $\mathcal{P}$. 
 \begin{defi}
If $s<0$, or \ $0<s<1$ a distribution $f \in \mathcal{S}'(\mathbb R^2)$ belongs to the homogeneous Besov space $\dot B^{s}_{\infty,\infty}$ if and only if, (in the case $0<s<1$, $\Delta_j f \in L^\infty$, for all $j\in \mathbb Z$)
$$\Vert f \Vert_{\dot B^{s}_{\infty,\infty}}=\sup_{j\in \mathbb Z} 2^{js}\Vert \Delta_j f \Vert_{\infty}< \infty,  $$
and $f \in \mathcal{S}'(\mathbb R^2)$ belongs to the homogeneous Besov space $\dot B^{0}_{\infty,1}$ if and only if
$$
\Vert f \Vert_{\dot B^{0}_{\infty,1}}=\sum_{j \in \mathbb Z } \Vert \Delta_j f \Vert_{\infty}< \infty.
$$
\end{defi}
 
\noindent A very useful tool when we work with Besov norms is the Bernstein inequalities.
 \begin{lemm} 
If $f \in \mathcal{S}'(\mathbb{R}^{2})$, then for all $(s, j) \in \mathbb{R} \times \mathbb{Z}$, and for all $1\leq p \leq q \leq \infty$ we have 
$$2^{js} \Vert \Delta_j f \Vert_{L^{p}} \lesssim \Vert \Lambda^{s} \Delta_{j} f \Vert_{p} \lesssim 2^{js} \Vert \Delta_j f \Vert_{L^{p}}, $$
$$\Vert \Lambda^{s} S_j f \Vert_{L^{p}} \lesssim 2^{js} \Vert S_j f \Vert_{L^{p}}, $$
\ \text{and} $$ \  \Vert \Delta_j f \Vert_q \lesssim2^{{j}(\frac{2}{p}-\frac{2}{q})}  \Vert \Delta_j f \Vert_{p}.$$
\end{lemm}

\section{The equation with truncated and regularized initial data.}
This section is devoted to the introduction of the $(SQG)_{R,\eps}$ equation which is the  critical (SQG) equation with truncated and regularized initial data.  As we said before, we put some integrability on $\theta$ by setting $\theta=\Lambda^{s} w$ and we consider the equation with respect to $w$ (which is more regular than $\theta$  by construction). We first truncate the initial data $w_{0}$ by considering    $w_{0} \chi_R$, where $\chi_R$  is a positive test function in $\mathcal{C}^\infty(\mathbb R^2)$ constructed as follows. Let $\chi \in \mathcal{D}(\mathbb R^{2})$ be a positive smooth function s.t $\chi(x)=1$ if $\vert x \vert \leq 1$, and 0 if $\vert x \vert \geq 2$. For $R>0$, we introduce the function  $\chi_R (x)\equiv \chi(x/R)$. Then we defined the truncated initial data $\theta_{0,R}$ by
$$
\theta_{0,R}=\Lambda^s w_{0,R}=\Lambda^s( w_{0}\chi_R).
$$

We will also need to regularize the initial data via a convolution with a standard mollifier $\rho_{\ep}$ because we need to have at least $H^{1}$ solutions. Namely, we consider a positive test function $\rho \in \mathcal{D}(\mathbb{R}^2)$, s.t. $supp (\rho) \subset [-1, 1]^2$ and the integral over $\mathbb R^2$ of $\rho$ is equal to one. Then we define $\rho_{\ep}(x)\equiv\ep^{-2}\rho (x\eps^{-1}).$ Therefore, we  consider the (SQG) equation associated with the following truncated and regularized initial data
$$
\theta_{0,R,\eps}=\Lambda^{s}(w_{0}\chi_{R}) * \rho_{\eps}.
$$
More precisely, we will focus on the following $(SQG)_{R,\ep}$ equation
\begin{equation}%\label{qg}
(SQG)_{R, \eps} :\ \\\left\{
\aligned
&\partial_{t}w_{{R,\ep}}=\left(\Lambda^{-s}\nabla \right)\cdot \left(\Lambda^{s}w_{{R,\ep}}\  \mathcal{R}^{\perp}{\Lambda^{s} w_{{R,\ep}}}\right) -\Lambda{w_{{R,\ep}}},
\\ \nonumber
& \nabla.\ \mathcal{R}^{\perp}{\Lambda^{s} w_{{R,\ep}}}=0, \hspace{4,9cm} 
\\ \nonumber
& \theta_{0,R,\eps}=\Lambda^{s}(w_{0}\chi_{R}) * \rho_{\eps}. 
\endaligned
\right.
\end{equation}
The condition on the truncated initial data is:
 $$\theta_{0,R}=\Lambda^s( w\chi_R) \in \Lambda^s(  H^{s}(\mathbb R^2)) \cap  L^{\infty}(\mathbb R^2) \subset L^{2}(\mathbb R^2) \cap L^{\infty}(\mathbb R^2).$$ Therefore, the truncated and regularized initial data satisfy 
$$
\displaystyle\theta_{0,R, \eps} \in \displaystyle \bigcap_{k\geq 0} H^{k} \subset H^{2} \subset  \dot B^0_{\infty,1}(\mathbb R^2).
$$ 
From the result of Resnick (see \cite{Res}), we can claim that there exists at least one solution $\theta_{R, \ep}$ which is a global weak Leray-Hopf solution, that is
$$
\theta_{R,\eps} \in L^{\infty}((0,T), L^{2}(\mathbb R^2)) \cap L^{2}((0,T),\dot H^{1/2}(\mathbb R^2)).
$$

\hspace{-0,5cm}Therefore, since $\theta_{R,\eps}=\Lambda^{s} w_{R,\eps},$ we also have the existence of a solution $w_{R,\ep}$ satisfies 

$$
w_{R,\eps} \in L^{\infty}((0,T),\dot H^{s}(\mathbb R^2)) \cap L^{2}((0,T),\dot H^{s+1/2}(\mathbb R^2)).
$$
Moreover, the result of Abidi and Hmidi \cite{Hm} (for data  in $\dot B^0_{\infty,1}(\mathbb R^2)$) provides us the existence of a global solution $\theta_{R,\ep}$ with the following regularity
$$
\theta_{R,\eps} \in{\mathcal{C}}(\RR_+;\,\dot B^0_{\infty,1})
\cap L^1_{\textnormal loc}(\RR_+;\,\dot B^1_{\infty,1}).
  $$
  The idea of the proof  is to obtain an energy estimate with respect to the parameters $R$ and $\ep$ for the solutions $w_{R,\ep}$. This energy inequality will provides us the desired compactness which allows us to pass to the limit in the $(SQG)_{R,\ep}$ equation. We will omit to treat the uniform bound for the initial data since it is already done in \cite{Laz}. Nevertheless, it is worth to recall the uniform estimates for the initial data as well as a useful bound for the Riesz transforms. The proof of those estimations can be found in \cite{Laz}.
  
  \newpage
\begin{lemm} \label{ine1} The truncated and regularized  initial data and the Riesz transforms  satisfy the following statements.  \\

$\bullet$ $ u_{R,\ep}$ is  bounded in $(L^{2}_{t}L^{2})_{uloc}$, moreover we have
$$
\Vert u_{R,\ep} \Vert_{(L^{2}_{t}L^{2})_{uloc}} \lesssim  \Vert w_{R, \eps} \Vert_{L^{2}_{t} \dot H^{s}_{uloc}}.
$$ 
$\hspace{0,5cm}\bullet$  If $w_{0,R,\eps} \in  H^{s}_{uloc}(\mathbb{R}^2)$ and $\theta_{0,R, \eps}=\Lambda^{s}w_{0,R, \eps} \in L^\infty(\mathbb{R}^2)$ then $$\sup_{R>1,\ep>0} \Vert w_{0,R,\eps}\Vert_{L^{\infty}(\mathbb{R}^2)}<\infty \ \  {\rm{and}} \ \ \sup_{R>1,\ep>0} \Vert w_{0,R,\eps}\Vert_{H^{s}_{uloc}(\mathbb{R}^2)}<\infty.$$
\end{lemm}

 We shall also use the $L^{p}$ maximum principle due to Resnick \cite{Res} and also C\'ordoba and C\'ordoba {\cite{CC}}. Namely, if $\theta$ is a smooth solution to the   $(SQG)_\alpha$ equation then all $p \in [2, +\infty)$ and for all $t \geq 0$

$$ 
\Vert \theta\Vert^{p}_{L^p}+2k\int^{t}_{0}\int \vert \Lambda^{\alpha/2}(\vert \theta\vert^{p/2}) \vert^2 \ dx \ ds \leq \Vert \theta_{0}\Vert^{p}_{L^p} \ \ {\rm{and }} \ \ \Vert \theta(t)\Vert_{L^\infty}\leq \Vert \theta_{0}\Vert_{L^\infty}.
$$

  \section{Global existence, the case $1/4\leq s \leq 1$.}
\noindent  In this section, we improve the following energy balance obtained  in \cite{Laz}, 
\begin{equation} \label{eq:equation1}
\partial_{t} A_{\phi}(w_{R,\ep}) \lesssim \Vert w_{R,\ep} \Vert^{2}_{  H^{s}_{uloc}(\mathbb{R}^2)} 
\end{equation}
where $w_{R,\ep}$ is a weak solution of $(SQG)_{R,\ep}$ and $A_{\phi} w$ is defined as in \ref{norm}. This inequality leads to a global existence result for all $1/2<s<1$ (see \cite{Laz}). By taking advantage of one term, we are able to prove a new energy inequality. This new energy inequality allows us to  extend the global existence result to the case $1/4\leq s \leq 1$. More precisely, we prove the following proposition

  \begin{prop} \label{tt}
Assume that $1/4 \leq s \leq 1$ and let $w_{R,\ep}$ be a weak solution of the equation $(SQG)_{R,\ep}$. Then, for all $\gamma>0$, for all $\phi \in B_{\phi_{0}}$,   we have the following energy inequality 
\begin{equation*}
\partial_{t} A_{\phi}(w_{R, \ep})  + \int \phi\theta_{R, \ep} \Lambda \theta_{R,\ep} \ dx \leq (C+\frac{1}{2\gamma} ) \ \Vert w_{R, \ep} \Vert^{2}_{ H^{s}_{uloc}(\mathbb{R}^2)} + \frac{\gamma}{2}  \Vert \phi w_{R,\ep} \Vert^{2}_{ \dot H^{s+1/2}},
\end{equation*}
where $C$ is a positive constant depending only on $\Vert \theta_{0} \Vert_{{L^\infty}(\mathbb{R}^2)}$ and $ \Vert w_{0} \Vert_{ H^{s}_{uloc}(\mathbb{R}^2)}.$ 
\end{prop}
\noindent{\bf Proof of  proposition \ref{tt}.}  Since 
$$
\Vert w_{R,\ep} \Vert^{2}_{ H^{s}_{uloc}}= \sup_{\phi \in B_{\phi_{0}}} \int \frac{ \vert \phi w_{R,\eps} \vert^2}{2} + \frac{ \vert \Lambda^{s}(\phi w_{R,\eps}) \vert^{2}}{2} \ dx,
$$
then, due to \ref{norm}, it is convenient to study the evolution of the quantity
$$
A_{\phi} w_{R,\eps} = \int \frac{ \vert \phi w_{R,\eps} \vert^2}{2} + \frac{ \vert \Lambda^{s}(\phi w_{R,\eps}) \vert^{2}}{2} \ dx.  
$$
A straighforward computation gives us the following equality 
\begin{align}
\label{uno}\partial_{t} A_{\phi}(w_{{R,\ep}}) +  \int \phi \Lambda^{s}w_{{R,\ep}} \Lambda^{s+1}w_{{R,\ep}} \  dx =&-\int  w_{{R,\ep}}\phi \Lambda^{-s}\nabla.(u_{{R,\ep}}\Lambda^{s}w_{{R,\ep}}) \ dx \\
&- \int \phi \Lambda^{s}w_{{R,\ep}} \nabla.(u_{{R,\ep}}\Lambda^{s}w_{{R,\ep}})   \ dx- \int w_{{R,\ep}}\phi   \Lambda{w_{{R,\ep}}} \ dx. \nonumber \end{align}

%$\partial_{t} A_{\phi}(w_{{R,\ep}}) +  \int \phi \Lambda^{s}w_{{R,\ep}} \Lambda^{s+1}w_{{R,\ep}} \  dx =-\int  w_{{R,\ep}}\phi \Lambda^{-s}\nabla.(u\Lambda^{s}w_{{R,\ep}}) \ dx \  $ \ \hspace{2cm} (1)\\
%\begin{flushright} $ - \int \phi \Lambda^{s}w_{{R,\ep}} \nabla.(u_{{R,\ep}}\Lambda^{s}w_{{R,\ep}})   \ dx- \int w_{{R,\ep}}\phi   \Lambda{w_{{R,\ep}}} \ dx$
%\end{flushright}
In \cite{Laz}, we proved that the last  two terms of the right-hand side of equality \eqref{uno}  are uniformly controlled by $C \Vert w_{R,\ep} \Vert^{2}_{\dot H^{s}_{uloc}}.$ The discussion about low and high oscillations (which, roughly speaking, corresponds respectively to small $s$ and big $s$) comes from the first term of the right hand side of equality \eqref{uno}, namely 
\begin{eqnarray*}
 -\int \Lambda^{-s}\nabla.(w_{{R,\ep}}\phi) u_{{R,\ep}}\Lambda^{s}w_{{R,\ep}} \ dx.  
\end{eqnarray*}
Since $\phi w_{R,\ep} \in L^{2}_{t}\dot H^{s+1/2}$ then $\Lambda^{-s}\nabla.(w_{{R,\ep}}\phi) \in L^{2}_{t} \dot H^{2s-1/2} $ we thus need the condition  $2s-1/2\geq0$ which implies $1/4 \leq s\leq 1$ (in some sense, integration by parts are forbidden since it would amount to control the $L^{\infty}$ norm of $u_{R,\ep}$ which is hopeless). Thus if $1/4 \leq s\leq 1$, we have by using the duality in homogeneous Sobolev  spaces  and by introducing a positive test function $\psi \in \mathcal{D}$, which is equal to 1 in a neighborhood of $supp (\phi)$  and 0 far enough from $supp (\phi)$, and then writting $\psi=\psi^{2}_1$  and  finally using Young's inequality, we obtain that for all $\gamma>0$
 \begin{eqnarray*}
 -\int \Lambda^{-s}\nabla.(w_{{R,\ep}}\phi) \psi  u_{{R,\ep}}\Lambda^{s}w_{{R,\ep}} \ dx  &\leq& \Vert  \Lambda^{-s}\nabla.(w_{{R,\ep}}\phi) \Vert_{\dot H^{2s-1/2}} \Vert  \psi u_{R,\ep} \theta_{R,\ep} \Vert_{\dot H^{1/2-2s}}  \\
 &\lesssim& \Vert \phi w_{R,\ep} \Vert_{\dot H^{s+1/2}}  \Vert \psi_1 u_{R,\ep} \Vert_{ L^{2}}  \Vert \theta_{0,{R,\ep}} \Vert_{L^\infty} \\
 &\lesssim& \Vert \phi w_{R,\ep} \Vert_{\dot H^{s+1/2}}  \Vert \psi_1w_{{R,\ep}} \Vert_{ H^{s}}  \Vert  \theta_{0,{R,\ep}} \Vert_{L^\infty} \\
 &\lesssim&\frac{\gamma}{2}  \Vert \phi w_{R,\ep} \Vert^{2}_{\dot H^{s+1/2}}+ \frac{1}{2\gamma}  \Vert  \psi_{1}w_{R,\ep} \Vert^{2}_{H^{s}} \\
  &\lesssim&\frac{\gamma}{2}  \Vert \phi w_{R,\ep} \Vert^{2}_{\dot H^{s+1/2}}+ \frac{1}{2\gamma} \underbrace{\sup_{\bar{\psi} \in B_{\psi_1}} \Vert  \bar{\psi}w_{R,\ep} \Vert^{2}_{ H^{s}}}_{=\Vert w_{R,\ep}\Vert^{2}_{H^{s}_{uloc}}}.
\end{eqnarray*}

Since the kernel of the operator $\Lambda^{-s}\nabla$ (that is $1/\vert x \vert^{3-s}$) is in $L^{1}(\mathbb R^2)$ far from the origin and since $u_{{R,\ep}}\Lambda^{s}w_{{R,\ep}} \in  L^{2}_{uloc}$ (because $u_{{R,\ep}} \in  L^{2}_{uloc}$ and futhermore it is controlled by $C \Vert w_{R,\ep} \Vert_{\dot H^{s}_{uloc}}$ and  by the maximum principle $\Lambda^{s}w_{{R,\ep}} \in L^\infty$), we find by using the Cauchy-Schwarz inequality
  \begin{eqnarray} \label{far}
 -\int w_{{R,\ep}}\phi \Lambda^{-s}\nabla.((1-\psi)  u_{{R,\ep}}\Lambda^{s}w_{{R,\ep}}) \ dx  &\lesssim& \Vert w_{{R,\ep}}\phi  \Vert_{L^{2}} \Vert u_{R,\ep} \Vert_{L^{2}_{uloc}}\lesssim \Vert w_{R,\ep} \Vert^{2}_{\dot H^{s}_{uloc}}.
  \end{eqnarray}
  
 \noindent Therefore, from  equality $\eqref{uno}$ we get the following inequality
    \begin{eqnarray*} \label{inn}
\partial_{t} A_{\phi}(w_{{R,\ep}}) +  \int \phi \Lambda^{s}w_{{R,\ep}} \Lambda^{s+1}w_{{R,\ep}} \  dx  \leq \frac{C\gamma}{2}  \Vert \phi w_{R,\ep} \Vert^{2}_{\dot H^{s+1/2}}+ (C+\frac{1}{2\gamma})  \Vert  w_{R,\ep} \Vert^{2}_{\dot H^{s}_{uloc}}.
  \end{eqnarray*}
  
 \noindent This completes the  proof of  Prop. \ref{tt}.
  
  \qed
  
 Now, we integrate the energy inequality of Prop. \ref{tt}  in time $s \in [0,T]$  
\begin{eqnarray*}
A_{\phi}(w_{{R,\ep}}(T))+ \int_{0}^{T} \int \phi \theta_{R,\ep} \Lambda \theta_{R,\ep} \ dx \ ds \leq A_{\phi}(w_{0,{R,\ep}}) &+&\frac{C\gamma}{2}  \int_{0}^{T}  \Vert \phi w_{R,\ep} \Vert^{2}_{\dot H^{s+1/2}} \ ds \\   && \hspace{-1cm} +  \       (C+\frac{1}{2\gamma} )  \int_{0}^{T} \Vert w_{R,\ep} \Vert^{2}_{ \dot H^{s}_{uloc}(\mathbb{R}^2)} \ ds.
\end{eqnarray*}
The last term of the left hand side provides us the $(L^{2} \dot H^{s+1/2})_{uloc}$ norm since we have the equality  
$$
\int_{0}^{T} \int \phi \theta_{R,\ep} \Lambda \theta_{R,\ep}\ dx \ ds=\int_{0}^{T} \int \Lambda^{1/2} \theta_{R,\ep} [\Lambda^{1/2}, \phi] \theta_{R,\ep} \ dx \ ds + \int_{0}^{T} \int \phi \vert \Lambda^{1/2} \theta_{R,\ep} \vert^{2} \ dx \ ds.
$$
Hence, we get
\begin{eqnarray*}
A_{\phi}(w_{{R,\ep}}(T)) &+& \int_{0}^{T} \int \phi \vert \Lambda^{1/2} \theta_{R,\ep} \vert^{2} \  dx \ ds \leq  A_{\phi}(w_{0,{R,\ep}}) \ +   \frac{C\gamma}{2}  \int  \Vert \phi w_{R,\ep} \Vert^{2}_{\dot H^{s+1/2}} \ ds        \\ &+& (C+\frac{1}{2\gamma} )  \int_{0}^{T} \Vert w_{R,\ep} \Vert^{2}_{ H^{s}_{uloc}(\mathbb{R}^2)} \ ds + \left\vert \int_{0}^{T} \int \Lambda^{1/2} \theta_{R,\ep} [\Lambda^{1/2}, \phi] \theta_{R,\ep} \ dx \ ds \right\vert.
\end{eqnarray*}
Then, we use the following estimate, for all $\nu>0$ and $\eta>0$ (see \cite{Laz} for the proof)
\begin{eqnarray*}
\left\vert \int_{0}^{T} \int \Lambda^{1/2} \theta_{R,\ep} [\Lambda^{1/2}, \phi] \theta_{R,\ep} \ dx \ ds \right\vert
& \lesssim & ( \frac{\nu}{2}+\frac{\eta}{2}) \int_{0}^{T} \Vert \psi \Lambda^{1/2} \theta_{R,\ep} \Vert^{2}_{L^2} \ ds+ (\frac{1}{2\nu}+\frac{1}{2\eta}) \int_{0}^{T} \Vert \theta_{R,\ep} \Vert^{2}_{L^{2}_{uloc}} \ ds. 
\end{eqnarray*}
where $\psi$ is an arbitrary positive test function which is equal to 1 in a neighborhood of the support of $\phi$. Finally, we obtain the following inequality
\begin{align} \label{ab}
A_{\phi}(w_{{R,\ep}}(T)) + \int_{0}^{T} \int \phi \vert \Lambda^{1/2} \theta_{R,\ep} \vert^{2} \ dx \ ds \leq & \ A_{\phi}(w_{0,R,\ep}) \ + C( \frac{\nu}{2}+ \frac{\eta}{2}) \int_{0}^{T}  \Vert \psi \Lambda^{1/2} \theta_{R,\ep} \Vert^{2}_{L^2} \ ds \nonumber \\ &+ \frac{C\gamma}{2}  \int  \Vert \phi w_{R,\ep} \Vert^{2}_{\dot H^{s+1/2}} \ ds   \\ 
& +  \left(C+\frac{1}{2\gamma}+\frac{1}{2\nu}+ \frac{1}{2\eta} \right)  \int_{0}^{T} \Vert w_{R,\ep} \Vert^{2}_{ H^{s}_{uloc}(\mathbb{R}^2)} \ ds.  \nonumber \\  & \hspace{-2cm}\ \leq A_{\phi}(w_{0,R,\ep}) \ + \left(C+\frac{1}{2\gamma}+\frac{1}{2\nu}+ \frac{1}{2\eta} \right)  \int_{0}^{T} \Vert w_{R,\ep} \Vert^{2}_{ H^{s}_{uloc}(\mathbb{R}^2)} \ ds \nonumber \\ &\hspace{-4cm}+ C( \frac{\nu}{2}+ \frac{\eta}{2}) \underbrace{\sup_{\bar{\psi} \in B_{\psi}} \int_{0}^{T}\Vert \bar{\psi} \Lambda^{1/2} \theta_{R,\ep} \Vert^{2}_{L^2} \ ds}_{=\Vert w_{R,\ep} \Vert^{2}_{(L^{2}\dot H^{s+1/2})_{uloc}}} \nonumber + \frac{C\gamma}{2}  \underbrace{\sup_{\phi \in B_{\phi_0}}\int_{0}^{T}  \Vert \phi w_{R,\ep} \Vert^{2}_{\dot H^{s+1/2}}}_{=\Vert w_{R,\ep} \Vert^{2}_{(L^{2}\dot H^{s+1/2})_{uloc}}}  \ ds .  \nonumber
\end{align}

\noindent Then, we take the supremum over all $\phi \in B_{\phi_0}$ and we choose $\eta$, $\nu$, $\gamma$ sufficiently small so that the norms $(L^{2}\dot H^{s+1/2})_{uloc}$ of the right hand side of (\ref{ab}) are absorbed by that of the left. 
\begin{remark}
It is important to note that  $\sup_{\phi \in B_{\phi_0}}  \int_{0}^{T} \int \phi \vert \Lambda^{1/2} \theta_{R,\ep} \vert^{2} \ dx \ ds = \Vert w_{R,\ep} \Vert^{2}_{(L^{2}\dot H^{s+1/2})_{uloc}}$
\end{remark}

\noindent We eventually get
\begin{align} \label{ro}
\hspace{1cm }\Vert w_{R,\ep} (T) \Vert^{2}_{ H^{s}_{uloc}(\mathbb{R}^2)} +C \Vert w_{R,\ep} \Vert^{2}_{ (L^{2}([0,T],\dot H^{s+1/2}(\mathbb{R}^2)))_{uloc}} \ dx  \  &\leq  \sup_{R>0,\ep>0} \Vert w_{0,R,\ep} \Vert^{2}_{ H^{s}_{uloc}(\mathbb{R}^2)} \\
&+ \ C  \int_{0}^{T} \Vert w_{R,\ep} \Vert^{2}_{  H^{s}_{uloc}(\mathbb{R}^2)} \ ds.   \nonumber
\end{align}
In particular, one obtain
 \begin{eqnarray*}
\Vert w_{{R,\ep}}(T) \Vert^{2}_{ H^{s}_{uloc}(\mathbb{R}^2)} &\leq& \sup_{R>0,\ep>0} \Vert w_{0,R,\ep} \Vert^{2}_{  H^{s}_{uloc}(\mathbb{R}^2)} +C  \int^{T}_{0} \Vert w_{R,\ep}(x,s) \Vert^{2}_{ H^{s}_{uloc}(\mathbb{R}^2)} \ ds. 
\end{eqnarray*}
 By  Gronwall's lemma we conclude that
\begin{equation*}
\hspace{1,8cm}  \Vert w_{{R,\ep}}\Vert^{2}_{L^{\infty}([0,T] ,  H^{s}_{uloc}(\mathbb{R}^2))} \lesssim e^{CT}. 
\end{equation*} 
On the other hand,  inequality \eqref{ro} gives  
\begin{equation*}
\Vert w_{R,\ep} \Vert^{2}_{ (L^{2} \dot H^{s+1/2})_{uloc}(\mathbb{R}^2)}   \lesssim \  \sup_{R>0,\ep>0}\Vert w_{0,R,\ep} \Vert^{2}_{ H^{s}_{uloc}(\mathbb{R}^2)} \ +C  \int_{0}^{T} \Vert w_{R,\ep} \Vert^{2}_{ \dot H^{s}_{uloc}(\mathbb{R}^2)} \ ds.
\end{equation*}
Thus,  $w_{R,\ep} \in  L^{2}([0,T],\dot H^{s+1/2}_{uloc}(\mathbb{R}^2))$  for all finite time $T>0$. For the limiting process when $R\rightarrow +\infty$ and $\ep \rightarrow 0$ we refer to the last section of \cite{Laz}. This completes the proof of the global existence result stated in Theorem \ref{tp}. The next section is devoted to the proof of the local existence.

\section{Local existence, the case $0<s<1/4$.}
\subsection{ Towards the energy inequality. } \label{rrr}

In the previous section (see Prop \ref{tt}), we have proved that, for all $1/4\leq s \leq 1$ and $\gamma>0$,
\begin{equation*}
\partial_{t} A_{\phi}(w_{R,\eps})  + \int \phi\theta_{R, \ep} \Lambda \theta_{R,\eps} \ dx \leq (C+\frac{1}{2\gamma} ) \ \Vert w_{R,\eps} \Vert^{2}_{  H^{s}_{uloc}(\mathbb{R}^2)} + \frac{\gamma}{2}  \Vert \phi w_{R,\eps} \Vert^{2}_{\dot H^{s+1/2}},
\end{equation*}
and this inequality leads to a global existence result.  These kind of estimates do not work for small values of $s$, that is, for $0<s <1/4$. Nevertheless, one can still obtain an energy inequality which provides a local existence result. As we have already recalled in the previous section, the main issue comes from the first term of the right hand side of  inequality $\eqref{uno}$. By integrating by parts, we see that this term can be rewritten  as
$$
-\int  w_{{R,\ep}}\phi \Lambda^{-s}\nabla.(u_{R,\ep}\Lambda^{s}w_{{R,\ep}}) \ dx = -\int  \Lambda^{-s}\nabla(w_{{R,\ep}}\phi) u_{R,\ep}\Lambda^{s}w_{{R,\ep}} \ dx.
$$
As before, we introduce a positive test function $\psi \in \mathcal{D}$  which is equal to 1 in a neighborhood of the support of $\phi$ and equal to 0 on supp$(\phi) +2$. Using H\"older's inequality, we get
$$-\int  \Lambda^{-s}\nabla(w_{{R,\ep}}\phi)  \psi u_{R,\ep}  \Lambda^{s}w_{{R,\ep}} \ dx \leq \Vert \phi w_{{R,\ep}} \Vert_{\dot H^{1-s}} \Vert \theta_{R,\eps} \Vert_{L^\infty} \Vert \psi u_{R,\eps} \Vert_{L^{2}}.$$
Unfortunately, this inequality fails for small values of $s$  (i.e. $0<s< 1/4$) since the right hand side is not controlled (note that we only have  $w_{R,\eps} \in L^{2} \dot H^{s+1/2}$). Basically, one would like to share the derivatives into the two terms so that both are controlled by a desired norm. In our case, if we share the derivatives (that is to say, if we integrate by parts) so that both terms are controlled i.e we give $s+1/2$ derivatives to $\phi w_{R,\ep}$ and $1/2-2s$ to $\psi u_{R,\ep}\Lambda^{s}w_{{R,\ep}}$, we obtain
$$-\int  w_{{R,\ep}}\phi \Lambda^{-s}\nabla \cdot (\psi u_{R,\eps}\Lambda^{s}w_{{R,\ep}}) \ dx= -\int  \Lambda^{s+\frac{1}{2}}(w_{{R,\ep}}\phi) \Lambda^{\frac{1}{2}-2s}\frac{\nabla}{\vert \nabla \vert}.(\psi u_{R,\ep}\Lambda^{s}w_{{R,\ep}}) \ dx.$$
Then, one would like to use the inequality (see for instance \cite{BCD}), where we set $\psi=\psi^2_1$
\begin{equation} \label{rep}
\Vert \psi u_{R,\ep} \theta_{R,\ep} \Vert_{\dot H^{\frac{1}{2}-2s}} \lesssim \Vert \psi_1 u_{R,\ep} \Vert_{L^\infty}\Vert \psi_1 \theta_{R,\ep} \Vert_{\dot H^{\frac{1}{2}-2s}}+\Vert \psi_1 \theta_{R,\ep} \Vert_{L^\infty}\Vert \psi_1 u_{R,\ep}  \Vert_{\dot H^{\frac{1}{2}-2s}}, 
\end{equation}
but this inequality clearly fails due to the lack of bound of  $\Vert \psi_1 u_{R,\ep} \Vert_{L^\infty}$. \\

To overcome the difficulty, we will show that one can still control the product  $\psi u_{R,\eps} \theta_{R,\eps}$ in a Sobolev space of positive regularity provided that we lose a little bit of regularity on $\psi u_{R,\eps}$. The idea consists in substituting the $L^\infty$ norm of $\phi u_{R,\eps}$ for the Besov norm $B^{-\delta}_{\infty, \infty}$ of $\phi u_{R,\eps}$, with  $\delta>0$. Indeed, we have the following lemma
\begin{lemm}
\label{producte} For all $\phi \in B_{\phi_{0}}$, there exists  $0<\delta<3/2$ such that $\phi u_{R,\eps} \in  B^{-\delta}_{\infty, \infty}$.
\end{lemm}
\noindent{\bf Proof.}   Since $\phi u_{R,\ep} \in H^{1/2}$, by using  Bernstein's inequality twice (in $2D$), we get $H^{1/2}= B^{1/2}_{2,2} \subset B^{-1/2}_{2, \infty} \subset B^{-3/2}_{\infty, \infty}$. Moreover, we have that  $\phi u_{R,\ep} \in B^{0}_{\infty, \infty}$. Finally,  by interpolation we conclude that $\phi u_{R, \ep} \in \dot B^{-\delta}_{\infty, \infty}$ for all $0<\delta<3/2$.\\ 
\qed 

\noindent Then, for $\sigma>0$ we write
\begin{eqnarray*}
-\int  \Lambda^{-s}\nabla(w_{{R,\ep}}\phi) \psi u_{R,\ep}\Lambda^{s}w_{{R,\ep}} \ dx&=&-\int w_{{R,\ep}}\phi \Lambda^{-s}\nabla.((1-\psi)  u_{{R,\ep}}\Lambda^{s}w_{{R,\ep}}) \ dx     \\  &-& \int  \Lambda^{-\sigma-s}\nabla(w_{{R,\ep}}\phi)  \Lambda^{\sigma} (\psi_1 u_{R,\ep} \psi_1 \Lambda^{s}w_{{R,\ep}}) \ dx , \\
&=& (I)+(II).
\end{eqnarray*}
Observe that $(I)$ can be estimated as in  step \ref{far} (see proof of  Prop. \ref{tt}) since in this step we just use the fact that the kernel $1/\vert x \vert^{3-s} \in L^{1}(\mathbb R^2)$ far from the origin (note that we are outside the support of $\phi$ thus far from the origin) and this holds for all $0<s<1$ and in particular it holds for $0<s\leq1/4$ therefore we still have the same estimates as \ref{far} and we get
\begin{equation} \label{3t}
(I) \lesssim \Vert w_{R,\ep} \Vert^{2}_{ H^{s}_{uloc}}
\end{equation}

\noindent For $(II)$, we use the following lemma to estimate the product  $\psi u_{R,\eps}\theta_{R,\eps}$ (recall that $\theta_{R,\eps}=\Lambda^{s}w_{{R,\ep}}$).

\begin{lemm} \label{pr}
There exist $\sigma>0$ and $\delta>0$ with $\sigma+\delta<1/2$ such that $\psi_1 \theta_{R,\eps} \in \dot H^{\sigma+\delta} \cap L^\infty$ and $\psi_1 u_{R,\eps} \in \dot H^{\sigma} \cap B^{-\delta}_{\infty, \infty}$. Moreover, we have  $\psi u_{R,\eps} \theta \in \dot H^{\sigma}$ and  the estimate
$$
\Vert \psi u_{R,\eps}\theta_{R,\eps} \Vert_{\dot H^{\sigma}} \leq \Vert \psi_1 u_{R,\eps} \Vert_{B^{-\delta}_{\infty, \infty}} \Vert \psi_1 \theta_{R,\eps} \Vert_{\dot H^{\sigma+\delta}} + \Vert \psi_1 \theta_{R,\eps} \Vert_{L^{\infty}} \Vert \psi_1 \theta_{R,\eps} \Vert_{\dot H^{\sigma}} 
$$
\end{lemm}

\noindent{\bf Proof.}  The first part is nothing but Lemma \ref{producte} and the regularity condition on $\theta_{R,\eps}$ and $u_{R,\eps}$. The inequality is a consequence of the Littlewood-Paley decomposition for functions of positive Sobolev regularity whose Fourier transform is supported on a ball. The proof is based on the following fact (see e.g. Theorem 4.1 p 34  \cite{PGLR}). If  $E$  is a shift invariant Banach space of distribution, then, for all $1\leq q \leq \infty$, and all $(s, \tau)$ such that $-s < \tau < 0$ and $s>0$, the pointwise multiplication is a bounded bilinear operator from $(B^{s}_{q, E} \cap B^{\tau}_{\infty, \infty}) \times (B^{s}_{q, E} \cap B^{\tau}_{\infty, \infty})$ to $B^{s+\tau}_{q, E}$. Where, for all $\kappa \in \mathbb R$, such that $\kappa_0<\kappa<\kappa_1$ and $s= (1-\eta) \kappa_0+\eta \kappa_1$ and for all $\eta \in (0,1)$, the space $B^{\kappa}_{q, E}$ is defined by the following interpolation space
$$
B^{\kappa}_{q, E}=[H^{\kappa_0}_{E}, H^{\kappa_1}_{E}]_{\eta, q}.
$$
Let us also recall that, for  all $\kappa\in \mathbb R $,  $H^{\kappa}_{E}$ is the space $(Id-\Delta)^{-\frac{\kappa}{2}}E$ endowed with the norm  $$\Vert f \Vert_{H^{\kappa}_{E}}=\Vert (Id-\Delta)^{\frac{\kappa}{2}} f \Vert_{E}.$$
Therefore, for all functions $(f,g) \in (B^{s}_{q, E} \cap B^{\tau}_{\infty, \infty})^{2}$ and for all $1\leq q \leq \infty,$  we have 
\begin{align}\label{est}\Vert f g \Vert_{B^{s+\tau}_{q, E}} \leq \Vert f \Vert_{B^{s}_{q, E}} \Vert g \Vert_{ B^{\tau}_{\infty, \infty}} + \Vert f \Vert_{ B^{\tau}_{\infty, \infty}} \Vert g \Vert_{B^{s}_{q, E}}. 
\end{align}
%\hspace{2,5cm} (5)$$
If we choose $E=L^{2}$ and $q=2$ and use the fact that $\Vert \theta_{R,\eps}\Vert_{ B^{-\tau}_{\infty, \infty}} \leq  \Vert \theta_{R,\eps} \Vert_{ L^{\infty}}$ we conclude the proof by applying inequality \eqref{est}  with   $s=\sigma+\delta$ and $\tau=-\delta$.

\qed

The next section is devoted to the proof of the energy inequality.

\subsection{The energy inequality}

In this section, we prove the following proposition
\begin{prop} \label{energy}Let $w_{R,\ep}$ be a solution of $(SQG)_{R,\ep}$  equation. Then, for all  $K \in (0,2)$
\begin{eqnarray}
 \hspace{1cm}\Vert w_{R,\ep}(T) \Vert^{2}_{H^{s}_{uloc}} \leq \Vert w_{0,R,\ep} \Vert^{2}_{H^{s}_{uloc}} +   C \int_{0}^{T} \left(\Vert w_{R, \ep} \Vert^{2}_{ H^{s}_{uloc}}+\Vert w_{R, \ep} \Vert^{2(\frac{3-K}{2-K})}_{{\dot H^s}_{uloc}}  +  C \right) \ ds 
 \end{eqnarray}
 \end{prop}
 \begin{remark}
The constant $C>0$ depends only on the  norms $\Vert \theta_{0} \Vert_{L^{\infty}}$ and $\Vert w_{0} \Vert_{ H^{s}_{uloc}}$.
 \end{remark}
  \noindent For the sake of clarity, we omit to write the parameter $\eps$ in the following. 

 \subsubsection{Proof of  Proposition \ref{energy}}
 \hspace{-0,6cm} Let us first choose $\sigma>0$  such that
$$
\frac{1}{2}-2s < \sigma < \frac{1}{2},
$$
\hspace{-0,0cm}and then choose $\delta>0$  such that 
$$
\delta+\sigma< \frac{1}{2}.
$$
These conditions and  Lemma \ref{pr} allow us to write
\begin{eqnarray*}
-\int  \Lambda^{-\sigma-s}\nabla(w_{{R}}\phi)  \Lambda^{\sigma} (\psi u_{R} \Lambda^{s}w_{{R}}) &\leq& \Vert w_{{R}}\phi \Vert_{\dot H^{1-s-\sigma}} \Vert \psi u_R \theta_R\Vert_{\dot H^{\sigma}} \\
&\lesssim&  \Vert w_{{R}}\phi \Vert_{\dot H^{1-s-\sigma}} \Vert \psi_1 u_R \Vert_{ B^{-\delta}_{\infty, \infty}} \Vert \psi_1 \theta_R \Vert_{H^{\sigma+\delta}} \\ && + \ \Vert w_{{R}}\phi \Vert_{\dot H^{1-s-\sigma}} \Vert \psi_1 \theta_R \Vert_{L^{\infty}} \Vert \psi_1 \theta_R \Vert_{\dot H^{\sigma}},
\end{eqnarray*}
 where  $\psi=\psi^{2}_1$. Then we write,
$$
\Vert \psi_1 u_R \Vert_{ B^{-\delta}_{\infty, \infty}} \leq \Vert u_R \Vert_{B^{-\delta}_{\infty, \infty}} \leq \Vert S_0 u_R \Vert_{L^\infty} + \Vert (Id-S_0) u_R \Vert_{ B^{-\delta}_{\infty, \infty}},
$$

\hspace{-0,6cm}where we have used that

$$\Vert S_0 u_R \Vert_{B^{-\delta}_{\infty, \infty}} \leq  \Vert S_0 u_R \Vert_{L^\infty}.$$

\hspace{-0,7cm} Then, we observe that
\begin{eqnarray*}
\Vert (S_0 u_R) (x) \Vert_{L^\infty} \leq \sum_{k \in \mathbb{Z}^2 \backslash \{0,0\}} \frac{1}{\vert k \vert^3} \left(\int_{k + x} \vert u_R (y) \vert^2 \ dy \right)^{1/2} &\leq& C \sup_{k \in \mathbb{Z}^2 \backslash \{0,0\}}\left(\int_{k + x} \vert u_R (y) \vert^2 \ dy \right)^{1/2} \\ &\lesssim& \Vert u_R \Vert_{L^{2}_{uloc}} \\ &\lesssim& \Vert w_R \Vert_{\dot H^{s}_{uloc}}, 
\end{eqnarray*}
where, in the last inequality, we have used the first point of Lemma \ref{ine1}. Therefore,
$$
\Vert S_0 u_R \Vert_{L^\infty} \leq \Vert u_R \Vert_{L^{2}_{uloc}} \leq \Vert w_R \Vert_{\dot H^{s}_{uloc}.}
$$
For high frequencies we have, using  Bernstein's inequality
$$
\Vert (Id-S_0) u_R \Vert_{ B^{-\delta}_{\infty, \infty}}  \leq \sum_{j>0} \Vert \Lambda^{-\delta} \mathcal{R}^{\perp} \Delta_j \theta_R \Vert_{L^\infty} \leq \sum_{j>0} 2^{-j\delta} \Vert \theta_R \Vert_{L^\infty} \leq C \Vert \theta_R \Vert_{L^\infty},
$$
and we obtain,
$$
\Vert \psi_1 u_R \Vert_{B^{-\delta}_{\infty, \infty}} \leq C \Vert w_R \Vert_{\dot H^{s}_{uloc}} + \Vert \theta_R \Vert_{L^\infty} \leq C \Vert w_R \Vert_{\dot H^{s}_{uloc}} + \Vert \theta_{0,R} \Vert_{L^\infty}.
$$

\hspace{-0,6cm}Thus, the inequality

\begin{eqnarray*}
-\int  \Lambda^{-\sigma-s}\nabla(w_{{R}}\phi)  \Lambda^{\sigma} (\psi u_{R} \Lambda^{s}w_{{R}}) &\lesssim&       \Vert w_{{R}}\phi \Vert_{\dot H^{1-s-\sigma}} \Vert \psi_1 \theta_R \Vert_{L^{\infty}} \Vert \psi_1 \theta_R \Vert_{\dot H^{\sigma}}  \\ && + \ \Vert w_{{R}}\phi \Vert_{\dot H^{1-s-\sigma}} \Vert \psi_1 u_R \Vert_{ B^{-\delta}_{\infty, \infty}} \Vert \psi_1 \theta_R \Vert_{\dot H^{\sigma+\delta}},
\end{eqnarray*}

\hspace{-0,6cm}becomes
\begin{align}
-\int  \Lambda^{-\sigma-s}\nabla(w_{{R}}\phi)  \Lambda^{\sigma} (\psi u_{R} \Lambda^{s}w_{{R}}) \lesssim& \ \ \Vert w_{R}\phi \Vert_{\dot H^{1-s-\sigma}}  \Vert \psi_1 \theta_R \Vert_{\dot H^{\sigma}} \Vert \psi_1 \theta_R \Vert_{L^{\infty}} \nonumber\\  & + \Vert w_{{R}}\phi \Vert_{\dot H^{1-s-\sigma}} \Vert \psi_1 \theta_R \Vert_{\dot H^{\sigma+\delta}} \Vert \theta_{0,R} \Vert_{L^{\infty}} \nonumber \\
 & + \Vert w_{{R}}\phi \Vert_{\dot H^{1-s-\sigma}} \Vert w_R \Vert_{\dot H^{s}_{uloc}}  \Vert \psi_1 \theta_R \Vert_{\dot H^{\sigma+\delta}}. \nonumber
\end{align}

Then, we use the fact that $0<s\leq1/4$ and that the constant  $\sigma>0$ is chosen so that $1/2-2s<\sigma$ therefore $s<1-s-\sigma < s +1/2$, and  $\delta +\sigma<1/2$ by definition. These facts allow us to  interpolate  and write that for all $(\kappa, \nu,\eta)\in (0,1)^3$
\begin{eqnarray*}
\partial_t A_\phi w_R  + \int \phi \theta_R \Lambda^s \theta_R \ dx &\leq& C \Vert w_R \Vert^{2}_{ H^{s}_{uloc}} \\ &+& C \Vert w_{{R}}\phi \Vert^{1-\kappa}_{\dot H^{s}}\Vert w_{{R}}\phi \Vert^{\kappa}_{\dot H^{s+1/2}}       \Vert \psi_1 \theta_R \Vert^{1-\eta}_{L^{2}} \Vert \psi_1 \theta_R \Vert^{\eta}_{\dot H^{1/2}}             \Vert \theta_{0,R} \Vert_{L^{\infty}}  \\
&+& C\Vert w_{{R}}\phi \Vert^{1-\kappa}_{\dot H^{s}}\Vert w_{{R}}\phi \Vert^{\kappa}_{\dot H^{s+1/2}}   \      \Vert \psi_1 \theta_R \Vert^{\nu}_{\dot H^{1/2}} \Vert \psi_1 \theta_R \Vert^{1-\nu}_{L^2} \ \Vert \theta_{0,R}    \Vert_{L^{\infty}} \\
 &+&  C\Vert w_{{R}}\phi \Vert^{1-\kappa}_{\dot H^{s}}\Vert w_{{R}}\phi \Vert^{\kappa}_{\dot H^{s+1/2}}      \Vert  w_R \Vert_{\dot H^{s}_{uloc}}         \Vert \psi_1 \theta_R \Vert^{\eta}_{\dot H^{1/2}}   \Vert \psi_1 \theta_R \Vert^{1-\eta}_{L^{2}}.        
 \end{eqnarray*}
 \begin{remark}
 Note that the first term of the right hand side i.e. $C \Vert w_R \Vert^{2}_{H^{s}_{uloc}}$ comes from the control of the other terms (see \ref{3t})
 \end{remark}
 \noindent Now, we take the supremum over all $\phi \in B_{\phi_{0}}$ and all $\bar{\psi_1} \in B_{\psi_{1}}$ only on some norms, namely, those which will give a $L^{\infty}_t \dot H^{s}_{uloc}$ norm of $w_R$ or equivalently a $L^{\infty}_t L^{2}_{uloc}$ norm of  $\theta_{R}$. We obtain
 \begin{align*}
\partial_t A_\phi w_R  + \int \phi \theta_R \Lambda^s \theta_R \ dx \leq& \ C \Vert w_R \Vert^{2}_{H^{s}_{uloc}} \\
&+ C \sup_{\substack{\phi \in B_{\phi_0}\\ \bar{\psi_1} \in B_{\psi_1}}} \left\{\Vert w_{{R}}\phi \Vert^{1-\kappa}_{\dot H^{s}}  \Vert \bar{\psi_1} \theta_R \Vert^{1-\eta}_{L^{2}}  \right\} \Vert w_{{R}} \phi \Vert^{\kappa}_{\dot H^{s+1/2}}      \Vert \psi_1 \theta_R \Vert^{\eta}_{\dot H^{1/2}}             \Vert \theta_{0,R} \Vert_{L^{\infty}}  \\
&+ C  \sup_{\substack{\phi \in B_{\phi_0}\\ \bar{\psi_1} \in B_{\psi_1}}} \left\{ \Vert w_{{R}}\phi \Vert^{1-\kappa}_{\dot H^{s}} \Vert \bar{\psi_1} \theta_R \Vert^{1-\nu}_{L^2} \right\}         \Vert w_{{R}}\phi \Vert^{\kappa}_{\dot H^{s+1/2}}   \      \Vert \psi_1 \theta_R \Vert^{\nu}_{\dot H^{1/2}}  \ \Vert \theta_{0,R}    \Vert_{L^{\infty}} \\
 &+  C \sup_{\substack{\phi \in B_{\phi_0}\\ \bar{\psi_1} \in B_{\psi_1}}}\left\{\Vert w_{{R}}\phi \Vert^{1-\kappa}_{\dot H^{s}}  \Vert \bar{\psi_1} \theta_R \Vert^{1-\eta}_{L^{2}} \right\}  \Vert  w_R \Vert_{\dot H^{s}_{uloc}}         \Vert w_{{R}}\phi \Vert^{\kappa}_{\dot H^{s+1/2}}              \Vert \psi_1 \theta_R \Vert^{\eta}_{\dot H^{1/2}}.
 \end{align*}

 \noindent Therefore, we obtain that
 
  \begin{eqnarray*}
\partial_t A_\phi w_R  + \int \phi \theta_R \Lambda^s \theta_R \ dx &\leq& C \Vert w_R \Vert^{2}_{ H^{s}_{uloc}} +  \Vert w_{{R}}\phi \Vert^{2-\kappa-\eta}_{\dot H^{s}}     \     \Vert w_{{R}} \phi \Vert^{\kappa}_{\dot H^{s+1/2}}      \Vert \psi_1 \theta_R \Vert^{\eta}_{\dot H^{1/2}}             \Vert \theta_{0,R} \Vert_{L^{\infty}}  \\
&+& C   \Vert w_{{R}} \Vert^{2-\kappa-\nu}_{\dot H^{s}_{uloc}}          \Vert w_{{R}}\phi \Vert^{\kappa}_{\dot H^{s+1/2}}   \      \Vert \psi_1 \theta_R \Vert^{\nu}_{\dot H^{1/2}}  \ \Vert \theta_{0,R}    \Vert_{L^{\infty}} \\
 &+&  C\Vert w_{{R}} \Vert^{3-\eta-\kappa}_{\dot H^{s}_{uloc}}          \Vert w_{{R}}\phi \Vert^{\kappa}_{\dot H^{s+1/2}}              \Vert \psi_1 \theta_R \Vert^{\eta}_{\dot H^{1/2}}.          
 \end{eqnarray*} 
  \begin{remark} \label{rk}
 It is important to note that the norms of the type  $ \Vert w_{{R}}\phi \Vert_{\dot H^{s+1/2}}$ or $\Vert \psi_1 \theta_R \Vert_{\dot H^{1/2}}$ are the same only after having taken the supremum over all $\phi \in B_{\phi_{0}}$ and $\bar{\psi_1} \in B_{\psi_{1}}$ respectively. In our case, we have to take those suprema as late as we can (at least after having integrated in time) so that we get the desired $(L^{2}_{t} \dot H^{s+1/2})_{uloc}$ norm of $w_R$ or equivalently the $(L^{2}_t \dot H^{1/2})_{uloc}$ norm of $\theta_R$.
  \end{remark}

 Let us focus on the three last terms appearing in  the right hand side of the previous inequality, that is to say:
 
 $$
T_1\equiv  \Vert w_{{R}}\phi \Vert^{2-\kappa-\eta}_{\dot H^{s}_{uloc}}     \     \Vert w_{{R}} \phi \Vert^{\kappa}_{\dot H^{s+1/2}}      \Vert \psi_1 \theta_R \Vert^{\eta}_{\dot H^{1/2}}             \Vert \theta_{0,R} \Vert_{L^{\infty}},  
$$

 $$
T_2\equiv  \Vert w_{{R}}\phi \Vert^{2-\kappa-\nu}_{\dot H^{s}_{uloc}}     \     \Vert w_{{R}} \phi \Vert^{\kappa}_{\dot H^{s+1/2}}      \Vert \psi_1 \theta_R \Vert^{\nu}_{\dot H^{1/2}}            
 \Vert \theta_{0,R} \Vert_{L^{\infty}}, 
$$
and
$$
T_3 \equiv \Vert w_{{R}} \Vert^{3-\eta-\kappa}_{\dot H^{s}_{uloc}}          \Vert w_{{R}}\phi \Vert^{\kappa}_{\dot H^{s+1/2}}              \Vert \psi_1 \theta_R \Vert^{\eta}_{\dot H^{1/2}}.
$$
The estimations of the terms $T_1$ and $T_2$ are the same. Let us begin with $T_1$ 
$$
T_1 = \Vert w_{{R}}\phi \Vert^{2-\kappa-\eta}_{\dot H^{s}_{uloc}}   \Vert \theta_{0,R} \Vert_{L^{\infty}} \  \Vert w_{{R}} \phi \Vert^{\kappa}_{\dot H^{s+1/2}}      \Vert \psi_1 \theta_R \Vert^{\nu}_{\dot H^{1/2}}.
$$
Using  Young's inequality (with the exponent $2/(\kappa+\eta) >1$) we get for all $\rho_1>0$

$$
T_1 \leq \rho_{1}^{\frac{2}{\kappa+\eta-2} }\left(\frac{2-\kappa-\eta}{2} \right)\Vert w_{{R}}\phi \Vert^{2}_{\dot H^{s}_{uloc}}   \Vert \theta_{0,R} \Vert^{\frac{2}{2-\kappa-\eta}}_{L^{\infty}} +  \rho_{1}^{    \frac{2}{\kappa+\eta}}\left(\frac{\kappa+\eta}{2} \right)   \Vert w_{{R}} \phi \Vert^{\frac{2\kappa}{\kappa+\eta}}_{\dot H^{s+1/2}}   \Vert \psi_1 \theta_R \Vert^{\frac{2\eta}{\kappa+\eta}}_{\dot H^{1/2}}.
$$
Since, by Young's inequality  (with the exponent $\frac{\kappa+\eta}{\gamma}>1$) we have
\begin{eqnarray*}
\Vert \phi w_R \Vert^{\frac{2\kappa}{\kappa+\eta}}_{{\dot H^{s+\frac{1}{2}}}}    \Vert  \psi_1\theta_R \Vert^{\frac{2\eta}{\kappa+\eta}}_{{\dot H^{1/2}}}    &\leq&  \frac{\kappa}{\kappa+\eta} \Vert w_R \phi \Vert^{2}_{\dot H^{s+1/2}}+\frac{\eta}{\kappa+\eta}  \Vert \psi_1 \theta_R \Vert^{2}_{\dot H^{1/2}},
 \end{eqnarray*}
and thus, we obtain
$$
T_1 \leq \rho_{1}^{\frac{2}{\kappa+\eta-2}}\left(\frac{2-\kappa-\eta}{2} \right)\Vert w_{{R}}\phi \Vert^{2}_{\dot H^{s}_{uloc}}   \Vert \theta_{0,R} \Vert^{\frac{2}{2-\kappa-\eta}}_{L^{\infty}} +    \frac{\kappa \rho_{1}^{\frac{2}{\kappa+\eta}}}{2} \Vert w_R \phi \Vert^{2}_{\dot H^{s+1/2}}+ \frac{\eta\rho_{1}^{\frac{2}{\kappa+\eta}}}{2}  \Vert \psi_1 \theta_R \Vert^{2}_{\dot H^{1/2}}, 
$$
replacing $\nu$ by $\eta$ we infer that
$$
T_2 \leq \rho_{1}^{\frac{2}{\kappa+\nu-2}}\left(\frac{2-\kappa-\nu}{2} \right)\Vert w_{{R}}\phi \Vert^{2}_{\dot H^{s}_{uloc}}   \Vert \theta_{0,R} \Vert^{\frac{2}{2-\kappa-\nu}}_{L^{\infty}} +    \frac{\kappa \rho_{1}^{\frac{2}{\kappa+\nu}}}{2} \Vert w_R \phi \Vert^{2}_{\dot H^{s+1/2}}+ \frac{\nu\rho_{1}^{\frac{2}{\kappa+\nu}}}{2}  \Vert \psi_1 \theta_R \Vert^{2}_{\dot H^{1/2}}.
$$
Therefore, for some $A>0$ we have for $i=1,2$
$$T_{i} \lesssim \ \rho^{-A}  \Vert w_{{R}}\phi \Vert^{2}_{\dot H^{s}_{uloc}} + \rho^{A} \left( \Vert w_R \phi \Vert^{2}_{\dot H^{s+1/2}}+   \Vert \psi_1 \theta_R \Vert^{2}_{\dot H^{1/2}} \right). $$

\noindent The term $T_3$ is the one which is responsible for the local existence. Indeed, using  Young's inequality (with the exponent $2/(\kappa+\eta)>1$), we obtain that for all $\rho_{3}>0$ 

\begin{eqnarray*}
T_3 &\leq&\rho_{3}^{\frac{2}{\kappa+\eta}} \left(\frac{\kappa+\eta}{2} \right) \Vert \phi w_R \Vert^{\frac{2\kappa}{\kappa+\eta}}_{{\dot H^{s+\frac{1}{2}}}}    \Vert  \psi_1\theta_R \Vert^{\frac{2\eta}{\kappa+\eta}}_{{\dot H^{1/2}}}+ \ \rho_{3}^{\frac{2}{\kappa+\eta-2}}  \left( \frac{2-\kappa-\eta}{2} \right)  \Vert w_R \Vert^{\frac{2(3-\eta-\kappa)}{2-\kappa-\eta}}_{{\dot H^s}_{uloc}}.
 \end{eqnarray*}
Once again,  Young's inequality (with the exponent $\frac{\kappa+\eta}{\gamma}>1$)  gives
\begin{eqnarray*}
\Vert \phi w_R \Vert^{\frac{2\kappa}{\kappa+\eta}}_{{\dot H^{s+\frac{1}{2}}}}    \Vert  \psi_1\theta_R \Vert^{\frac{2\eta}{\kappa+\eta}}_{{\dot H^{1/2}}}    &\leq&  \frac{\kappa}{\kappa+\eta} \Vert w_R \phi \Vert^{2}_{\dot H^{s+1/2}}+\frac{\eta}{\kappa+\eta}  \Vert \psi_1 \theta_R \Vert^{2}_{\dot H^{1/2}}.
 \end{eqnarray*}
Therefore,
$$
T_3  \leq  \frac{\kappa \rho_{3}^{\frac{2}{\kappa+\eta}}}{2} \Vert w_R \phi \Vert^{2}_{\dot H^{s+1/2}}+ \frac{\eta \rho_{3}^{\frac{2}{\kappa+\eta}}}{2}  \Vert \psi_1 \theta_R \Vert^{2}_{\dot H^{1/2}}  + \rho_{3}^{\frac{2}{\kappa+\eta-2}} \left( \frac{2-\kappa-\eta}{2} \right) \Vert w_R \Vert^{\frac{2(3-\kappa-\eta)}{2-\kappa-\eta}}_{{\dot H^s}_{uloc}}.
$$
Hence, for some constant $B>0$ we have
$$
T_3  \lesssim \  \rho_{3}^{B} \Vert w_R \phi \Vert^{2}_{\dot H^{s+1/2}}+  \rho_{3}^{B}  \Vert \psi_1 \theta_R \Vert^{2}_{\dot H^{1/2}}  + \rho_{3}^{-B} \Vert w_R \Vert^{\frac{2(3-\kappa-\eta)}{2-\kappa-\eta}}_{{\dot H^s}_{uloc}},
$$

 \noindent Putting all those terms together leads us to the following inequality
 \begin{eqnarray*}
\partial_t A_\phi w_R  + \int \phi \theta_R \Lambda^s \theta_R \ dx &\leq& \left( C + \rho^{-A}  \right) \Vert w_R \Vert^{2}_{H^{s}_{uloc}} +  \rho_{3}^{-B} \Vert w_R \Vert^{\frac{2(3-\kappa-\eta)}{2-\kappa-\eta}}_{{\dot H^s}_{uloc}}  \\ &+& ( \rho^{A} + \rho_{3}^{B} ) \Vert w_R \phi \Vert^{2}_{\dot H^{s+1/2}}+  ( \rho^{A} + \rho_{3}^{B} ) \Vert \psi_1 \theta_R \Vert^{2}_{\dot H^{1/2}}.
   \end{eqnarray*}

 \hspace{-0,6cm}We integrate in time $s \in [0,T]$, and we conclude that
  \begin{align} \label{inegg}
 A_\phi w_R(x,T) + \int_{0}^{T} \int \theta_R \Lambda^s \theta_R \ dx  \ ds &\leq  A_\phi w_{0,R} (x,T) + \left( C + \rho^{-A}  \right)  \int_{0}^{T} \Vert w_R \Vert^{2}_{H^{s}_{uloc}} \ ds  \\ &+  ( \rho^{A} + \rho_{3}^{B} ) \left ( \int_{0}^{T} \Vert w_R \phi \Vert^{2}_{\dot H^{s+1/2}} \ ds +  \int_{0}^{T}  \Vert \psi_1 \theta_R \Vert^{2}_{\dot H^{1/2}} \ ds \right) \nonumber \\
 &+ \rho_{3}^{-B} \int_{0}^{T} \Vert w_R \Vert^{\frac{2(3-\gamma-\eta)}{2-\gamma-\eta}}_{{\dot H^s}_{uloc}} \ ds.  \nonumber
   \end{align}
As before, we write the last integral of the left hand side in terms of the commutator as  
 \begin{equation}
\int_{0}^{T} \int \phi \theta_R \Lambda \theta_R \ dx \ ds=\int_{0}^{T} \int \Lambda^{1/2} \theta_R [\Lambda^{1/2}, \phi] \theta_R \ dx \ ds + \int_{0}^{T} \int \phi \vert \Lambda^{1/2} \theta_R \vert^{2} \ dx \ ds. 
\end{equation}
The last integral of the right hand side provides a $(L^{2}_t \dot H^{s+1/2})_{uloc}$ norm of $w_R$ which allows us to absorbe the other $(L^{2}_t \dot H^{s+1/2})_{uloc}$ norms of $w_R$, or equivalently, the $(L^{2}_t \dot H^{1/2})_{uloc}$ norms of $\theta_R$ appearing in the right hand side of the inequality (this equality of norms is verified only after having taken the supremum over the test functions, see  Remark \ref{rk}). The reminder term of the equality (\ref{inegg}) is controlled, since we have seen that for all $\mu_1 >0$ and $\rho_4>0$, 
 \begin{eqnarray*}
\left\vert \int_{0}^{T} \int \Lambda^{1/2} \theta_{R} [\Lambda^{1/2}, \phi] \theta_{R} \ dx \ ds \right\vert
& \lesssim & ( \frac{\mu_1}{2}+\frac{\rho_4}{2}) \int_{0}^{T} \Vert \psi \Lambda^{1/2} \theta_{R} \Vert^{2}_{L^2} \ ds \\ && + \  (\frac{1}{2\mu_1}+\frac{1}{2\rho_4}) \int_{0}^{T} \Vert \theta_{R} \Vert^{2}_{L^{2}_{uloc}} \ ds. 
\end{eqnarray*}
Therefore,  inequality $\eqref{inegg}$ becomes
\begin{align}
 A_\phi w_R(x,T) + \int_{0}^{T} \int \phi \vert \Lambda^{1/2} \theta_R \vert^{2} & dx ds   \lesssim   A_\phi w_{0,R} (x,T) + \left( C + \rho^{-A}  \right)  \int_{0}^{T} \Vert w_R \Vert^{2}_{ H^{s}_{uloc}} \ ds   \\ &+  ( \rho^{A} + \rho_{3}^{B} ) \left ( \int_{0}^{T} \Vert w_R \phi \Vert^{2}_{\dot H^{s+1/2}} \ ds +  \int_{0}^{T}  \Vert \psi_1 \theta_R \Vert^{2}_{\dot H^{1/2}} \ ds \right) \nonumber \\
 &+ ( \frac{\mu_1}{2}+\frac{\rho_4}{2}) \int_{0}^{T} \Vert \psi \Lambda^{1/2} \theta_{R} \Vert^{2}_{L^2} \ ds+ (\frac{1}{2\mu_1}+\frac{1}{2\rho_4}) \int_{0}^{T} \Vert \theta_{R} \Vert^{2}_{L^{2}_{uloc}} \ ds  \nonumber
 \\ & +     \rho_{3}^{-B} \int_{0}^{T} \Vert w_R \Vert^{\frac{2(3-\gamma-\eta)}{2-\gamma-\eta}}_{{\dot H^s}_{uloc}} \ ds.  \nonumber
 \end{align}

 \hspace{-0,6cm}Now we take the supremum over $\phi \in B_{\phi_{0}}$, $\bar{\psi_1} \in B_{\psi_1}$ and  $\bar{\psi} \in B_{\psi}$ in the last inequality. Therefore, we obtain
 \begin{eqnarray*}
 \Vert w_R(T) \Vert^{2}_{ H^{s}_{uloc}} + \Vert w_R\Vert^{2}_{{(L^{2}\dot H^{s+1/2})}_{uloc}} &\lesssim& \Vert w_{0,R} \Vert^{2}_{ H^{s}_{uloc}} +  (C+ \rho^{-A}+ \frac{1}{2\mu_1}+\frac{1}{2\rho_4}) \int_{0}^{T} \Vert w_R \Vert^{2}_{H^{s}_{uloc}} \ ds    \\ &+&   ( \rho^{A} + \rho_{3}^{B} +  \frac{\mu_1}{2}+\frac{\rho_4}{2} )   \int_{0}^{T}  \Vert  w_R \Vert^{2}_{\dot H^{s+1/2}_{uloc}} \ ds  \\  &+&\rho_{3}^{-B} \int_{0}^{T} \Vert w_R \Vert^{\frac{2(3-\gamma-\eta)}{2-\gamma-\eta}}_{{\dot H^s}_{uloc}} \ ds.  
 \end{eqnarray*}
We choose $\rho$, $\rho_3$, $\rho_4$  and $\mu_1$ sufficiently small  so that the  ${{(L^{2}\dot H^{s+1/2})}_{uloc}}$ norms  appearing  in the right hand side of the above  inequality are absorbed by that of the left, we get
\begin{eqnarray} \label{ine}
 \Vert w_R(T) \Vert^{2}_{ H^{s}_{uloc}} + \Vert w_R\Vert^{2}_{{(L^{2}\dot H^{s+1/2})}_{uloc}} & \lesssim \Vert w_{0,R} \Vert^{2}_{H^{s}_{uloc}} +  C \int_{0}^{T} \Vert w_R \Vert^{2}_{ H^{s}_{uloc}} \ ds   \\ & +    \int_{0}^{T} \Vert w_R \Vert^{\frac{2(3-\gamma-\eta)}{2-\gamma-\eta}}_{{\dot H^s}_{uloc}} \ ds.  \nonumber  
 \end{eqnarray}
 Thus, for all $K\equiv \gamma+\eta \in (0,2)$ we conclude that
 \begin{eqnarray} \label{fin}
 \Vert w_R(T) \Vert^{2}_{ H^{s}_{uloc}} \leq \Vert w_{0,R} \Vert^{2}_{H^{s}_{uloc}} +   C \int_{0}^{T} \left(\Vert w_R \Vert^{2}_{ H^{s}_{uloc}}+\Vert w_R \Vert^{\frac{2(3-K)}{2-K}}_{{\dot H^s}_{uloc}} \right)   \ ds.
 \end{eqnarray}
 In particular, we have
  \begin{eqnarray*} 
 \Vert w_R(T) \Vert^{2}_{ H^{s}_{uloc}} \leq \Vert w_{0,R} \Vert^{2}_{H^{s}_{uloc}} +   C \int_{0}^{T} \left(\Vert w_R \Vert^{2}_{ H^{s}_{uloc}}+\Vert w_R \Vert^{\frac{2(3-K)}{2-K}}_{{\dot H^s}_{uloc}} \right)   \ ds.
 \end{eqnarray*}

  \noindent This completes the proof of Proposition \ref{energy}.
 
 \qed

\noindent Since inequality (\ref{fin}) is true for all $K \in (0,2)$, it is still verified e.g. for $K=1$  and we get
$$
 \hspace{1cm}\Vert w_R(T) \Vert^{2}_{ H^{s}_{uloc}} \leq \Vert w_{0,R} \Vert^{2}_{ H^{s}_{uloc}} +    \int_{0}^{T} \left(\Vert w_R \Vert^{2}_{ H^{s}_{uloc}}+\Vert w_R \Vert^{4}_{{\dot H^s}_{uloc}}  +  C\right) \ ds. 
$$

\noindent From the previous inequality, we can write
$$
\partial_t \Vert w_R(t) \Vert^{2}_{ H^{s}_{uloc}} \leq    C (\Vert \theta_0 \Vert_{\infty}) \Vert w_R \Vert^{2}_{H^{s}_{uloc}}+\Vert w_R \Vert^{4}_{{\dot H^s}_{uloc}}  +  C. 
$$
Let us set
$$
\Omega (t)\equiv\Vert w_R (t) \Vert^{2}_{H^{s}_{uloc}} +1.
$$
Then,
$$
\partial_t \Omega (t)=\partial_t \Vert w_R(t) \Vert^{2}_{ H^{s}_{uloc}} \leq C (\Vert \theta_0 \Vert_{\infty})(\Vert w_R \Vert^{2}_{H^{s}_{uloc}} +1)^{2}=C (\Vert \theta_0 \Vert_{\infty})(\Omega (t))^{2}.
$$
By integrating, we obtain 
$$
\frac{1}{\Omega (0)}-\frac{1}{\Omega (T)} \leq C (\Vert \theta_0 \Vert_{\infty}) T,
$$
therefore, we get
\begin{eqnarray} \label{br}
\Vert w_R \Vert^{2}_{L^{\infty} H^{s}_{uloc}} \leq \left(\frac{1}{1+\Vert w_{0,R} \Vert^{2}_{ H^{s}_{uloc}}}-C (\Vert \theta_0 \Vert_{\infty})T\right)^{-1}
\end{eqnarray}
Thus, the solutions exist in $L^{\infty} H^{s}_{uloc}$  as long as $$T< T^{\star} \equiv\frac{C (\Vert \theta_0 \Vert_{\infty})}{ 1+\Vert w_{0,R} \Vert^{2}_{ H^{s}_{uloc}}}$$

\noindent Now, we have to show that $w_R \in (L^{2}\dot H^{s+1/2})_{uloc}$. In inequality \eqref{ine}, we have seen that 
\begin{eqnarray*}
 \Vert w_R\Vert^{2}_{{(L^{2}\dot H^{s+1/2})}_{uloc}} &\leq& \Vert w_{0,R} \Vert^{2}_{H^{s}_{uloc}} +  C \int_{0}^{T} \Vert w_R \Vert^{2}_{H^{s}_{uloc}} \ ds   \\ &&+ \  C \int_{0}^{T} \Vert w_R \Vert^{4}_{{\dot H^s}_{uloc}} \ ds + CT.
\end{eqnarray*}
Therefore inequality \eqref{br} allows us to conclude that the $(L^{2}\dot H^{s+1/2})_{uloc}$ norm of  $w_R$ is controlled for all $T<C (\Vert \theta_0 \Vert_{\infty})  (1+\Vert w_{0,R} \Vert^{2}_{H^{s}_{uloc}})^{-1}$. \\

\noindent Concerning the passage to the weak limit with respect to the parameters $R$ and $\ep$, we refer to the last section of \cite{Laz} since we have obtained the same uniform bounds.

 \vskip0.2cm\noindent{\bf Acknowledgment}:  The author is grateful to Prof. Marco Cannone and Prof. Pierre-Gilles Lemari{\'e}-Rieusset for encouragement and fruitful discussions.

 \vspace{1cm}

\begin{tabular}{ll}
{\small Universit\'e Paris-Est, Marne-la-Vall\'ee} & {\small  Instituto de Ciencias Matem\'aticas}\\
{\small LAMA, UMR 8050 CNRS } & {\small Consejo Superior de Investigaciones Cient\'ificas   }\\
{\small }5, boulevard Descartes, Cit\'e Descartes     & {\small  C/ Nicol\'{a}s Cabrera, 13-15, Campus Cantoblanco UAM }\\
{\small }77454 Marne-la-Vall\'ee, France  & {\small  28049 Madrid, Spain}\\
{\small Email: omar.lazar@univ-mlv.fr} & {\small  omar.lazar@icmat.es }
\end{tabular}

\end{document}